\documentclass[11pt]{article}
\usepackage[english]{babel}
\usepackage{amsmath}
\usepackage{amsthm}
\usepackage{amsfonts}
\usepackage{amssymb}
\usepackage{enumerate}
\usepackage{xspace}
\usepackage{euscript}
\usepackage{graphicx}
\usepackage{amscd}
\usepackage{epsfig}
\usepackage{tabularx}
\usepackage{url}
\usepackage{multimedia}

\newcommand{\be}{\begin{eqnarray}}
\newcommand{\ben}{\begin{eqnarray*}}
\newcommand{\en}{\end{eqnarray}}
\newcommand{\enn}{\end{eqnarray*}}

\newcommand{\G}{\Gamma}

\title{A spectral projection method for transmission eigenvalues}

\author{Fang Zeng \thanks{Institute of Computing and Data Sciences,
College of Mathematics and Statistics,
           Chongqing University, Chongqing 401331, China ({\tt fzeng@cqu.edu.cn}).}
\and Jiguang Sun \thanks{Department of
Mathematical Sciences, Michigan Technological University, Houghton, MI 49931 U.S.A. ({\tt jiguangs@mtu.edu}).}
\and Liwe Xu \thanks{Institute of Computing and Data Sciences,
College of Mathematics and Statistics,
           Chongqing University, Chongqing 401331, China ({\tt xul@cqu.edu.cn}).}}
\date{}

\begin{document}
\maketitle
\baselineskip 3.0ex
\begin{abstract}
In this paper, we consider a nonlinear integral eigenvalue problem, which is a reformulation of
the transmission eigenvalue problem arising in the inverse scattering theory.
The boundary element method is employed for discretization, which leads to a generalized matrix eigenvalue problem.
We propose a novel method based on the spectral projection.
The method probes a given region on the complex plane using contour integrals and decides if the region contains eigenvalue(s) or not.
It is particularly suitable to test if zero is an eigenvalue of the generalized eigenvalue problem, which in turn implies that
the associated wavenumber is a transmission eigenvalue.
Effectiveness and efficiency of the new method are demonstrated by numerical examples.
\end{abstract}

\section{Introduction}
We consider a non-linear non-selfadjoint transmission eigenvalue problem, which arises in the inverse scattering theory \cite{ColtonKress2013, CakoniEtal2010IP}.
Since 2010, the problem has attracted quite some attention of numerical mathematicians \cite{ColtonMonkSun2010IP, Sun2011SIAMNA,
JiSunTurner2012ACMTOM, AnShen2013JSC, SunXu2013IP, Kleefeld2013IP, CakoniMonkSun2014CMAM, LiEtal2014JSC, HuangHuangLin2015SIAMSC}.
The first numerical treatment by Colton, Monk, and Sun appeared in
\cite{ColtonMonkSun2010IP}, 
where three finite element methods
were proposed. A mixed method was developed by Ji, Sun, and Turner in \cite{JiSunTurner2012ACMTOM}.
An and Shen \cite{AnShen2013JSC} 
proposed an efficient spectral-element based
numerical method for transmission 
eigenvalues of two-dimensional, radially-stratified media. 
The first method supported by a  
rigorous convergence analysis was introduced by Sun in
\cite{Sun2011SIAMNA}.
Recently, Cakoni et.al. \cite{CakoniMonkSun2014CMAM} 
reformulated the problem and proved convergence
(based on Osborn's compact operator theory 
\cite{Osborn1975MC}) of a mixed finite 
element method. Li et.al. 
\cite{LiEtal2014JSC}
developed a finite element method based 
on writing the TE as a quadratic eigenvalue problem.
Other methods 
\cite{GintidesPallikarakis2013IP, JiSun2013JCP,
JiSunXie2014JSC, GintidesPallikarakis2013IP, YanHanBi2016} 
have been proposed recently.

Despite significant effort to develop various numerical methods for the transmission eigenvalue problem, 
computation of both real and complex eigenvalues remains difficult due to the fact that the numerical discretization usually end up with large sparse
generalized non-Hermitian eigenvalue problems, which are very challenging in numerical linear algebra. Traditional methods such as shift and invert Arnoldi are handicapped by the lack of a priori spectrum information.

In this paper, we adopt an integral formulation for the transmission eigenvalue problem. Using boundary element method,
the integral equations are discretized and a generalized eigenvalue problem of dense matrices is obtained. The matrices are significantly smaller than those from finite element methods.
If zero is a generalized eigenvalue, the corresponding wavenumber is a transmission eigenvalue. We propose a probing method
based on the spectral projection using contour integrals. The closed contour is chosen to be a small circle centered at the origin and a numerical quadrature is used
to compute the spectral projection of a random vector. The norm of the projected vector is used as an indicator of whether zero is an eigenvalue or not.

Integral based
methods \cite{Goedecker1999RMP, SakuraiSugiura2003CAM,
Polizzi2009PRB, Beyn2012LAA}
for eigenvalue computation, having their roots in the classical spectral perturbation theory (see, e.g., \cite{Kato1966}),
become popular in many areas, e.g., electronic structure calculation. These
methods are based on eigenprojections
using contour integrals of
the resolvent
\cite{AustinKravajaTrefethen2014SIAMNA}. Randomly chosen functions are projected to the generalized eigenspace
corresponding to the eigenvalues inside a closed contour, which leads to a relative small finite dimension eigenvalue problem.
For recently developments along this line, we refer the readers to \cite{KramerEtal2013JCAM, TangPolizziSIAMMAA, YingChanYeung2014, Yin2015}

For most existing integral based methods, estimation on
the locations, number of eigenvalues and dimensions of eigenspace are critical for their successes.
The proposed method is related to the methods developed in \cite{Kleefeld2013IP} and \cite{HuangEtal2015}.
The rest of the paper is arranged as follows.
In Section 2, we introduce the transmission eigenvalue problem and rewrite it using integral operators.
In Section 3, we present the probing method based on contour integrals.
We present numerical results in Section 4. Discussion and future works are contained in Section 5. 

\section{The transmission eigenvalue problem}
Let $D \subset \mathbb R^2$ be an open bounded domain with $C^2$ boundary $\Gamma := \partial D$.
The transmission eigenvalue problem is to find $k \in \mathbb C$ such that there exist non-trivial solutions $w$ and $v$ satisfying
\begin{subequations}\label{ATE}
\begin{align}
\label{ATEw}\Delta w+k^2 n w&=0, &\text{in } D,\\[1mm]
\label{ATEv}\Delta v+k^2 v&=0, &\text{in } D,\\[1mm]
\label{ATEbcD}w -v &= 0, &\text{on } \Gamma, \\[1mm]
\label{ATEbcN} \frac{\partial w}{\partial \nu} - \frac{\partial v}{\partial \nu}&= 0, &\text{on } \Gamma,
\end{align}
\end{subequations}
where $\nu$ is the unit outward normal to $\Gamma$.
The wavenumber $k$'s for which the transmission
eigenvalue problem has non-trivial solutions
are called transmission eigenvalues.
Here $n$ is the index of refraction, which is assumed to be a constant greater than $1$ in this paper.
Note that, for the integral formulation to be used, the index of refraction needs to be constant (see, e.g., \cite{HsiaoEtal2011CAM}).


In the following, we describe an integral formulation of the transmission eigenvalue problem following \cite{CossoniereHaddar2013JIE} (see also \cite{Kleefeld2013IP}).
Let $\Phi_k$ be the Green's function given by
\[
\Phi_k(x, y) = \frac{i}{4} H_0^{(1)}(k|x-y|),
\]
where $H_0^{(1)}$ is the Hankel function of the first kind of order $0$. The single and double layer potentials are defined as
\begin{eqnarray*}
&&(SL_k\phi)(x) = \int_{\partial \Omega} \Phi_k(x,y) \phi(x) \,ds(y), \\
&&(DL_k \phi)(x) = \int_{\partial \Omega}\frac{\partial \Phi_k}{\partial \nu(y)}(x,y) \phi(x) \,ds(y),
\end{eqnarray*}
where $\phi$ is the density function.

Let $(v, w)\in H^1(D) \times H^1(D)$ be a solution to \eqref{ATE}. Denote by $k_1 = \sqrt{n} k$ and set
\begin{eqnarray*}
&&\alpha := \frac{\partial v} {\partial \nu} \Big|_{\Gamma} =  \frac{\partial w} {\partial \nu} \Big|_{\Gamma}  \in H^{-1/2}(\Gamma), \\
&&\beta:= v|_{\Gamma} = w|_{\Gamma} \in H^{1/2}(\Gamma).
\end{eqnarray*}
Then $v$ and $w$ has the following integral representation
\begin{subequations}\label{IRvw}
\begin{align}
\label{IRv}&v = SL_k \alpha - DL_k \beta, &\text{in } D,\\[1mm]
\label{IRw}&w = SL_{k_1} \alpha - DL_{k_1}\beta, &\text{in } D.
\end{align}
\end{subequations}

Let $u:=w-v$. Then $u|_{\Gamma} = 0$ and $\frac{\partial u}{\partial \nu} |_{\Gamma} = 0$.
The boundary conditions of \eqref{ATE} imply that the transmission eigenvalues are $k$'s such that
\begin{equation}\label{Zk}
Z(k)\begin{pmatrix} \alpha \\ \beta \end{pmatrix} = 0,
\end{equation}
where
\[
Z(k) = \begin{pmatrix}  S_{k_1} - S_k & - K_{k_1} + K_k \\ - K_{k_1}'+K_k' & T_{k_1} - T_{k} \end{pmatrix}
\]
and the potentials $S_k, K_k, K_k', T_k$ are given by
\begin{subequations}\label{SkKkKkTk}
\begin{align}
\label{Sk}&(S_k\phi)(x) = \int_{\Gamma} \Phi_k(x,y) \phi(y) ds(y), \\[1mm]
\label{Kk}&(K_k \psi)(x) = \int_{\Gamma} \frac{\partial \Phi_k}{\partial \nu(y)} (x,y) \phi(y) ds(y), \\[1mm]
\label{Sk}&(K_k' \phi)(x) = \int_{\Gamma} \frac{\partial \Phi_k}{\partial \nu(x)} (x,y) \phi(y) ds(y),\\[1mm]
\label{Kk}&(T_k \psi)(x) = \frac{\partial}{\partial \nu(x)} \int_{\Gamma} \frac{\partial \Phi_k}{\partial \nu(y)} (x,y) \phi(y) ds(y).
\end{align}
\end{subequations}

It is shown in \cite{CossoniereHaddar2013JIE} that
\[
Z(k):= H^{-3/2}(\Gamma) \times H^{-1/2}(\Gamma) \to H^{3/2}(\Gamma) \times H^{1/2}(\Gamma)
\]
is of Fredholm type with index zero and analytic on $\mathbb C \setminus \mathbb R^{-}$.

From \eqref{Zk}, $k$ is a transmission eigenvalue if zero is an eigenvalue of $Z(k)$. Unfortunately, $Z(k)$ is compact.
The eigenvalues of $Z(k)$ accumulate at zero, which makes it impossible to distinguish zero and other eigenvalues numerically.
The workaround proposed in \cite{Cossoniere2011} is to consider a generalized eigenvalue problem
\begin{equation}\label{AkXBkX}
Z(k)\begin{pmatrix} \alpha \\ \beta \end{pmatrix} = \lambda B(k) \begin{pmatrix} \alpha \\ \beta \end{pmatrix},
\end{equation}
where $B(k) = Z(ik)$.
Since there does not exist purely imaginary transmission eigenvalues \cite{ColtonMonkSun2010IP}, the accumulation
point is shifted to $-1$. Then $0$ becomes isolated.

Now we describe a boundary element discretization of the potentials and refer the readers to \cite{NISTHandbook, SauterSchwab2011} for more details.
One discretizes the boundary $\G$ into element segments. 
Suppose the computational boundary $\G$ is discretized into $N$ segments $\G_1,\G_2,...,\G_N$ by nodes $x_1,x_2,...,x_N$ and 
$\tilde{\G}=\cup_{i=1}^N\G_i$. Let $\{\psi_j\},j=1,2,...,N$, be piecewise constant basis functions and $\{\varphi_j\},j=1,2,...,N,$ be
piecewise linear basis functions. We seek an approximate solution $\alpha_h$ and $\beta_h$ in the form
\ben
\alpha_h=\sum_{j=1}^N{\alpha_j\psi_j},\quad \beta_h=\sum_{j=1}^N{\beta_j\varphi_j}.
\enn
We arrive at a linear system
\ben
(V_{k,h}-V_{k_1,h})\Vec{\alpha}+(-K_{k,h}+K_{k_1,h})\Vec{\beta} &=& 0,\\
(K'_{k,h}-K'_{k_1,h})\Vec{\alpha}+(W_{k,h}-W_{k_1,h})\Vec{\beta} &=& 0,
\enn
where $\Vec{\alpha}=(\alpha_1,...,\alpha_N)^T$, $\Vec{\beta}=(\beta_1,...,\beta_N)^T$, and $V_{k,h},K_{k,h},K'_{k,h},W_{k,h}$ are matrices with entries
\ben
V_{k,h}(i,j) &=& \int_{\tilde\G}{(S_k\psi_j)\psi_i ds},\\
K_{k,h}(i,j) &=& \int_{\tilde\G}{(K_k\varphi_j)\psi_i ds},\\
K'_{k,h}(i,j) &=& \int_{\tilde\G}{(K'_k\psi_j)\varphi_i ds},\\
W_{k,h}(i,j) &=& \int_{\tilde\G}{(T_k\varphi_j)\varphi_i ds}.
\enn
In the above matrices, we can use series expansions of the first kind Hankel function as
\ben
H_0^{(1)}(x) &=& \sum_{m=0}^{\infty}{\frac{(-1)^m}{(m!)^2} \left(\frac{x}{2}\right)^{2m}}+\frac{2i}{\pi}\sum_{m=0}^{\infty} {\frac{(-1)^m}{(m!)^2}\left(\frac{x}{2}\right)^{2m}}\left(\ln\frac{x}{2}+c_e\right)\\
&-& \frac{2i}{\pi}\sum_{m=0}^{\infty}{\frac{(-1)^m}{(m!)^2} \left(\frac{x}{2}\right)^{2m}}\left(1+\frac{1}{2}+\frac{1}{m}\right) ,
\enn
where $c_e$ is the Euler constant. Thus,
\ben
H_0^{(1)}(k|x-y|)=\sum_{m=0}^{\infty} \left(C_5(m)+C_6(m)\ln{\frac{k}{2}}\right)k^{2m}|x-y|^{2m} +C_6(m)\ln{|x-y|}k^{2m}|x-y|^{2m},
\enn
where
\ben
C_5(m)&=&\frac{(-1)^m}{2^{2m}(m!)^2} \left[1+\frac{2c_ei}{\pi}- \frac{2i}{\pi}\left(1+\frac{1}{2}+\frac{1}{m}\right)\right],\\
C_6(m)&=&\frac{(-1)^mi}{2^{2m-1}(m!)^2\pi}.
\enn
We also need the following integrals which can be computed exactly.
\ben
Int_7(m) &=& \int_{-1}^1\int_{-1}^1 (\xi_1-\xi_2)^{2m}d\xi_2d\xi_1 \\
&=& \frac{2^{2m+2}}{(2m+1)(m+1)},
\enn
\ben
Int_8(m) &=& \int_{-1}^1\int_{-1}^1 (\xi_1-\xi_2)^{2m}\ln|\xi_1-\xi_2|d\xi_2d\xi_1 \\
&=& \frac{2^{2m+2}\ln2}{(2m+1)(m+1)}-\frac{(4m+3)2^{2m+3}}{(2m+1)^2(2m+2)^2},
\enn
\ben
Int_9(m) &=& \int_{-1}^1\int_{-1}^1 (\xi_1-\xi_2)^{2m}\xi_1\xi_2d\xi_2d\xi_1 \\
&=& \sum_{l=0}^{2m}\frac{(-1)^lC_{2m}^l}{(l+2)(2m+2-l)}[1-(-1)^l]^2,
\enn
and
\ben
Int_{10}(m) &=& \int_{-1}^1\int_{-1}^1 (\xi_1-\xi_2)^{2m}\xi_1\xi_2\ln|\xi_1-\xi_2|d\xi_2d\xi_1 \\
&=& \frac{-m2^{2m+2}\ln2}{(2m+1)(m+1)(m+2)}+ \frac{1}{(2m+1)(m+1)} \left[\frac{2^{2m+3}}{2m+3}-\frac{2^{2m+2}}{(m+2)^2}-\frac{2^{2m+1}}{m+1}\right]\\
&+& \frac{1}{2(m+1)^2(2m+1)^2} \sum_{l=0}^{2m+1}C_{2m+1}^l \left[\frac{(2m+1)^2}{l+2}(1-(-1)^l)- \frac{4m+3}{l+3}(1-(-1)^{l+1})\right].
\enn

Now we consider
\ben
V_{k,h}(i,j) &=& \int_{\tilde\G}{(V_k\psi_j)\psi_i ds}\\
&=& \int_{\tilde\G}\int_{\tilde\G}\Phi_k(x,y)\psi_j(y)\psi_i(x) ds_yds_x\\
&=& \int_{\G_i}\int_{\G_j}\Phi_k(x,y)\psi_j(y)\psi_i(x) ds_yds_x.
\enn
The integral over $\G_i\times \G_j$ can be calculated as
\ben
\int_{\G_i}\int_{\G_j}\Phi_k(x,y)\psi_j(y)\psi_i(x) ds_yds_x &=& \frac{i}{4}\int_{\G_i}\int_{\G_j}H_0^{(1)}(k|x-y|)\psi_j(y)\psi_i(x) ds_yds_x\\
&=& \frac{iL_iL_j}{16}\int_{-1}^1\int_{-1}^1H_0^{(1)}(k|x(\xi_1)-y(\xi_2)|) d\xi_2d\xi_1,
\enn
where
\ben
x(\xi_1)=x_i+\frac{1+\xi_1}{2}(x_{i+1}-x_i),\\
y(\xi_2)=x_j+\frac{1+\xi_2}{2}(x_{j+1}-x_j).
\enn
When $i\ne j$, it can be calculated by Gaussian quadrature rule. When $i=j$, we have
\ben
&\quad& \frac{iL_i^2}{16}\int_{-1}^1\int_{-1}^1H_0^{(1)}(k|x(\xi_1)-y(\xi_2)|) d\xi_2d\xi_1 \\
&=& \frac{iL_i^2}{16}\sum_{m=0}^{\infty}\frac{k^{2m}L_i^{2m}}{2^{2m}} \left(C_5(m)+C_6(m)\ln\frac{kL^i}{4}\right) \int_{-1}^1\int_{-1}^1 (\xi_1-\xi_2)^{2m} d\xi_2d\xi_1 \\
&& \qquad + \frac{iL_i^2}{16}\sum_{m=0}^{\infty}\frac{k^{2m}L_i^{2m}}{2^{2m}} C_6(m) \int_{-1}^1\int_{-1}^1 (\xi_1-\xi_2)^{2m}\ln|\xi_1-\xi_2| d\xi_2d\xi_1\\
&=& \sum_{m=0}^{\infty}\frac{ik^{2m}L_i^{2m+2}}{2^{2m+4}} \left[\left(C_5(m)+C_6(m)\ln\frac{kL^i}{4}\right) Int_7(m)+C_6(m)Int_8(m)\right].
\enn
The following regularization formulation is needed to discretize the hyper-singular boundary integral operator 
\begin{equation}\label{Wkbetax}
W_k \beta(x)=-\frac{d}{ds_x}V_k(\frac{d\beta}{ds})(x)-k^2\nu_x\cdot V_k(\beta\nu)(x).
\end{equation}
We refer the readers to \cite{Hsiao2011} for details of the discretization.

The above boundary element method leads to the following generalized eigenvalue problem
\begin{equation}\label{AxlambdaBx}
A{\bf x} = \lambda B{\bf x},
\end{equation}
where $A, B \in \mathbb C^{n \times n}$,
$\lambda \in \mathbb C$ is a scalar, and
${\bf x} \in \mathbb C^{n}$.

To compute transmission eigenvalues, the following method is proposed in \cite{Cossoniere2011}.
A searching interval for wavenumbers is discretized.
For each $k$,  the boundary integral operators $Z(k)$ and $Z(ik)$ are discretized to obtain \eqref{AxlambdaBx}. Then all eigenvalues
$\lambda_i(k)$ of \eqref{AxlambdaBx} are computed and arranged such that
\[
0 \le |\lambda_1(k)| \le |\lambda_2(k)| \le \ldots
\]
If $k$ is a transmission eigenvalue, $|\lambda_1|$ is very close to $0$ numerically. If one plots the inverse of $|\lambda_1(k)|$ against $k$, the
transmission eigenvalues are located at spikes.

\section{The probing method}
The method in \cite{Cossoniere2011} only uses the smallest eigenvalue. 
Hence it is not necessary to compute all eigenvalues of \eqref{AkXBkX}.
In fact, there is no need to know the exact value of $\lambda_1$. The only thing we need is that, if $k$ is a transmission eigenvalue, the generalized
eigenvalue problem \eqref{AkXBkX} has an isolated eigenvalue close to $0$.
This motivates us to propose a probing method to test if $0$ is an generalized eigenvalue of \eqref{AkXBkX}. The method does not compute
the actual eigenvalue and only solves a couple of linear systems. The workload is reduced significantly in two dimension and even more in three dimension.

We start to recall some basic results from spectral theory of compact operators \cite{Kato1966}. 
Let $T: \mathcal{X} \to \mathcal{X} $ be a compact operator on a complex Hilbert space $\mathcal{X} $. The resolvent set of $T$ is defined as
\begin{equation}\label{rhoT}
\rho(T)=\{ z \in \mathbb C: (z-T)^{-1} \text{ exists as a bounded operator on } \mathcal{X} \}.
\end{equation}
For any $z \in \rho(T)$, the resolvent operator of $T$ is defined as
\begin{equation}\label{resolventT}
R_z(T) = (z-T)^{-1}.
\end{equation}
The spectrum of $T$ is $\sigma(T)=\mathbb C \setminus \rho(T)$.
We denote the null space of an operator $A$ by $N(A)$. Let $\alpha$ be such that
\[
N\left((\lambda-T)^\alpha\right) = N\left((\lambda-T)^{\alpha+1}\right).
\]
Then $m=\dim N\left((\lambda-T)^\alpha \right)$ is called the algebraic multiplicity of $\lambda$.
The vectors in $N\left((\lambda-T)^\alpha\right)$ are called generalized eigenvectors
of $T$ corresponding to $\lambda$. Geometric multiplicity of $\lambda$ is defined as $\dim N(\lambda - T)$.

Let $\gamma$ be a simple closed curve on the complex plane $\mathbb C$ lying in $\rho(T)$, which contains $m$ eigenvalues, counting multiplicity, of $T$: $\lambda_i, i = 1, \ldots, m$.
We set
\[
P= \frac{1}{2\pi i} \int_\gamma R_z(T) dz.
\]
It is well-known that $P$ is a projection from $ \mathcal{X}$ onto the space of generalized
eigenfunctions ${\bf u}_i, i=1, \ldots, m$
associated with $\lambda_i, i=1, \ldots,m$ \cite{Kato1966}.

Let ${\bf f} \in \mathcal{X}$ be randomly chosen. If there are no eigenvalues inside $\gamma$, we have that $P{\bf  f} = {\bf 0}$.
Therefore, $P{\bf f}$ can be used to decide if a region contains eigenvalues of $T$ or not.

For the generalized matrix eigenvalue problem \eqref{AxlambdaBx},
the resolvent is
\begin{equation}\label{eigresolvent}
R_z(A, B) = (zB-A)^{-1}
\end{equation}
for $z$ in the resolvent set of the matrix pencil $(A, B)$.
The projection onto the generalized
eigenspace corresponding to eigenvalues
enclosed by $\gamma$
is given by \begin{equation}\label{ABpro}
P_k(A, B) = \frac{1}{2\pi i} \int_\gamma (zB-A)^{-1} dz .
\end{equation}
We write $P_k$ to emphasize that $P$ depends on the wavenumber $k$.

The approximation of $P_k{\bf f}$ is computed by suitable quadrature rules
\begin{equation}\label{XLXf}
	P_k{\bf f} = \dfrac{1}{2 \pi i} \int_{\gamma} R_z(A, B){\bf f} {d} z
	\approx \dfrac{1}{2 \pi i} \sum_{j=1}^W \omega_j R_{z_j}(A, B) {\bf f} =  \dfrac{1}{2 \pi i} \sum_{j=1}^W \omega_j {\bf x}_j,
\end{equation}
where $w_j$ are weights and $z_j$ are quadrature points. Here ${\bf x}_j$'s are the solutions of the following linear systems
\begin{equation}\label{linearsys}
(z_jB - A){\bf x}_j = {\bf f}, \quad j = 1, \ldots, W.
\end{equation}
Similar to the continuous case, if there are no eigenvalues inside $\gamma$,
then $P_k=0$ and  thus $P_k{\bf f} = {\bf 0}$ for all
${\bf f} \in \mathbb C^n$.  Similar to \cite{HuangEtal2015},
we project the random vector twice for a better result, i.e., we compute $P^2_{k} {\bf f}$.

For a fixed wavenumber $k$, the algorithm of the probing method is as follows.
\begin{itemize}
\item[]{\bf Input:}  a small circle $\gamma$ center at the origin with radius $r \ll 1$ and a random ${\bf f}$
\item[]{\bf Output:}  0 - k is not a transmission eigenvalue; 1 - k is a transmission eigenvalue
\item[1.] Compute $P_k^2{\bf f}$ by \eqref{XLXf};
\item[2.] Decide if $\gamma$ contains an eigenvalue:
	\begin{itemize}
		\item No. output 0.
		\item Yes. output 1.
	\end{itemize}
\end{itemize}

\section{Numerical Examples}
We start with an interval $(a, b)$ of wavenumbers and uniformly divide it into $K$ subintervals. At each wavenumber
\[
k_j = a+jh, \quad j = 0, 1, \ldots, K, \, h = \frac{b-a}{K},
\]
we employ the boundary element method to discretize the potentials.
We choose $N=32$ and end up with a generalized eigenvalue problem \eqref{AxlambdaBx} with $64 \times 64$ matrices $A$ and $B$.
To test whether $0$ is a generalized eigenvalue of \eqref{AxlambdaBx}, we choose $\gamma$ to be a circle of radius $1/100$.
Then we use $16$ uniformly distributed quadrature points on $\gamma$ and evaluate the eigenprojection \eqref{XLXf}. 
If at a wavenumber $k_j$, the projection is approximately $1$, then $k_j$ is a transmission
eigenvalue. For the actual computation, we use a threshold value $\sigma = 1/2$ to decide if $k_j$ is a transmission eigenvalue or not, i.e.,
$k_j$ is a transmission eigenvalue if $\|P_{k_j}^2 {\bf f} \|/\|P_{k_j} {\bf f}\|\ge \sigma$ and not otherwise.

Let $D$ be a disk with radius $1/2$. The index of refraction is $n=16$.
In this case, the exact transmission eigenvalues are known \cite{ColtonMonkSun2010IP}.
They are $k$'s such that
\begin{equation}
\label{y0}J_1(k/2)J_0(2k) - 4 J_0(k/2)J_1(2k) = 0
\end{equation}
and
\begin{equation}
\label{ym}J_{m-1}(k/2) J_m(2k) - 4J_m(k/2) J_{m-1}(2k)=0
\end{equation} 
for  $m=1, 2, \ldots$.
The actual values are given in Table \eqref{TableITEDisk}.
\begin{center}
\begin{table}
\caption{TEs of a disk with radius $r=1/2$ and index of refraction $n=16$.}
\label{TableITEDisk}
\begin{center}
\begin{tabular}{|c|c|c|c|}
\hline
        {$m=0$}&{1.9880}&{3.7594}&{6.5810}\\
\hline
          {$m=1$}&{2.6129}&{4.2954}&{5.9875}\\
\hline
         {$m=2$}&{3.2240}&{4.9462}&{6.6083}\\
\hline
\end{tabular}
\end{center}
\end{table}
\end{center}

We choose the interval to be $(1.5, 3.5)$ and uniformly divide it into $2000$ subintervals. At each $k_j$ we
compute the projection \eqref{XLXf} twice. The probing method finds three eigenvalues in $(1.5, 3.5)$
\[
k_1 = 1.988, \quad k_2 = 2.614, \quad k_3 = 3.228,
\]
which approximate the exact eigenvalues (the first column of Table  \eqref{TableITEDisk}) accurately.
Note that the continuous finite element method in \cite{ColtonMonkSun2010IP} computes 
\[
k_1=2.0301, \quad k_2=2.6937, \quad k_3= 3.3744,
\]
on a triangular mesh with mesh size $\approx 0.1$. The method proposed in this paper is more accurate.
However, we would like to remark that the methodology of the finite element method in \cite{ColtonMonkSun2010IP} is totally different.

\begin{figure}
\begin{center}
{ \scalebox{0.6} {\includegraphics{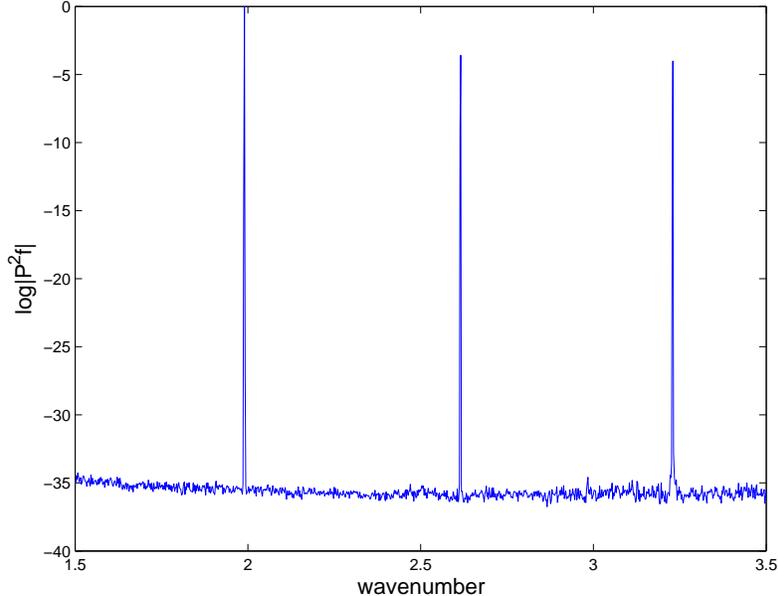}}}
\caption{The plot of $\log|P^2 {\bf f}|$ against the wavenumber $k$ for $n=16$.}
\label{n16}
\end{center}
\end{figure}

We also plot the log of $|P^2{\bf f}|$ against the wavenumber $k$ in Fig. \ref{n16}.
The method is robust since the eigenvalues can be easily identified.

\begin{figure}
\begin{center}
{ \scalebox{0.6} {\includegraphics{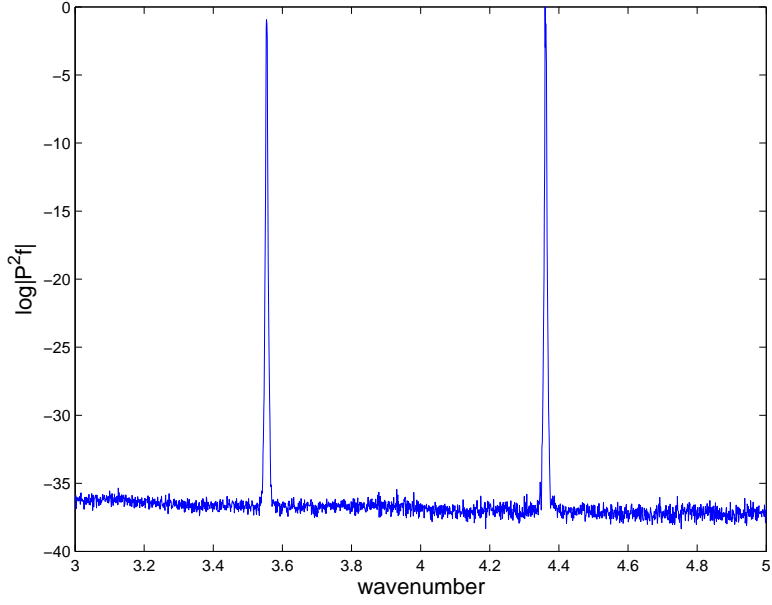}}}
\caption{The plot of $\log|P^2 {\bf f}|$ against the wavenumber for $n=9$.}
\label{n9}
\end{center}
\end{figure}

We repeat the experiment by choosing $n=9$ and $(a, b) = (3, 5)$. The rest parameters keep the same. The following eigenvalues are obtained
\[
k_1 = 3.554, \quad k_2 = 4.360.
\]
The log of $|P^2{\bf f}|$ against the wavenumber $k$ is shown in Fig.~\ref{n9}.

Finally, we compare the proposed method with the method in \cite{Cossoniere2011}. We take $n=16$ and compute for 2000 wavenumbers.
The CPU time in second is shown in Table~\ref{Time}. Note that all the computation is done using Matlab R2014a on a MacBook Pro with a 3 GHz Intel Core i7 and 16 GB memory.
We can see that the proposed method saves more time if the size of the generalized
eigenvalue problem is larger. We expect that it has a greater advantage for three dimension problems since the size of
the matrices are much larger than two dimension cases.
\begin{center}
\begin{table}
\caption{Comparison. The first column is the size of the matrix problem. The second column is the time used by the proposed method in second.
The second column is the time used by the method given in \cite{Cossoniere2011}. The fourth column is the ratio.}
\label{Time}
\begin{center}
\begin{tabular}{|l|r|r|r|}
\hline
          {size}&{probing method}&{method in \cite{Cossoniere2011}}& ratio\\
\hline
        {$64 \times 64$ }&{1.741340}&{5.742839}&3.30\\
\hline
         {$128 \times 128$}&{5.653961}&{31.152448}&5.51\\
\hline
         {$256 \times 256$}&{25.524530}&{224.435704}&8.79\\
\hline
         {$512 \times 512$}&{130.099433}&{1822.545973}&14.01\\
\hline
\end{tabular}
\end{center}
\end{table}
\end{center}

\begin{figure}
\begin{center}
\begin{tabular}{cc}
\scalebox{0.35}{\includegraphics{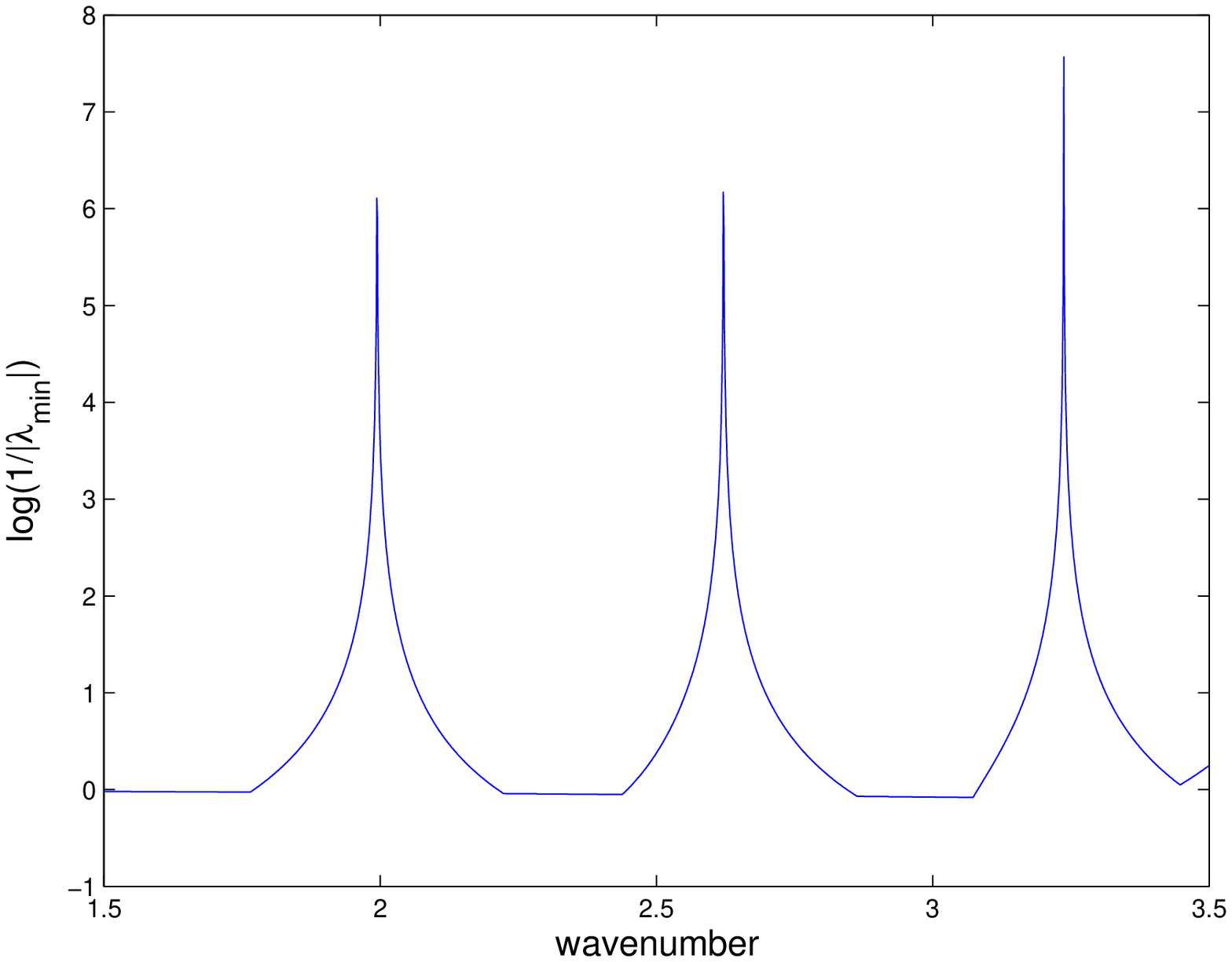}}&
\scalebox{0.35}{\includegraphics{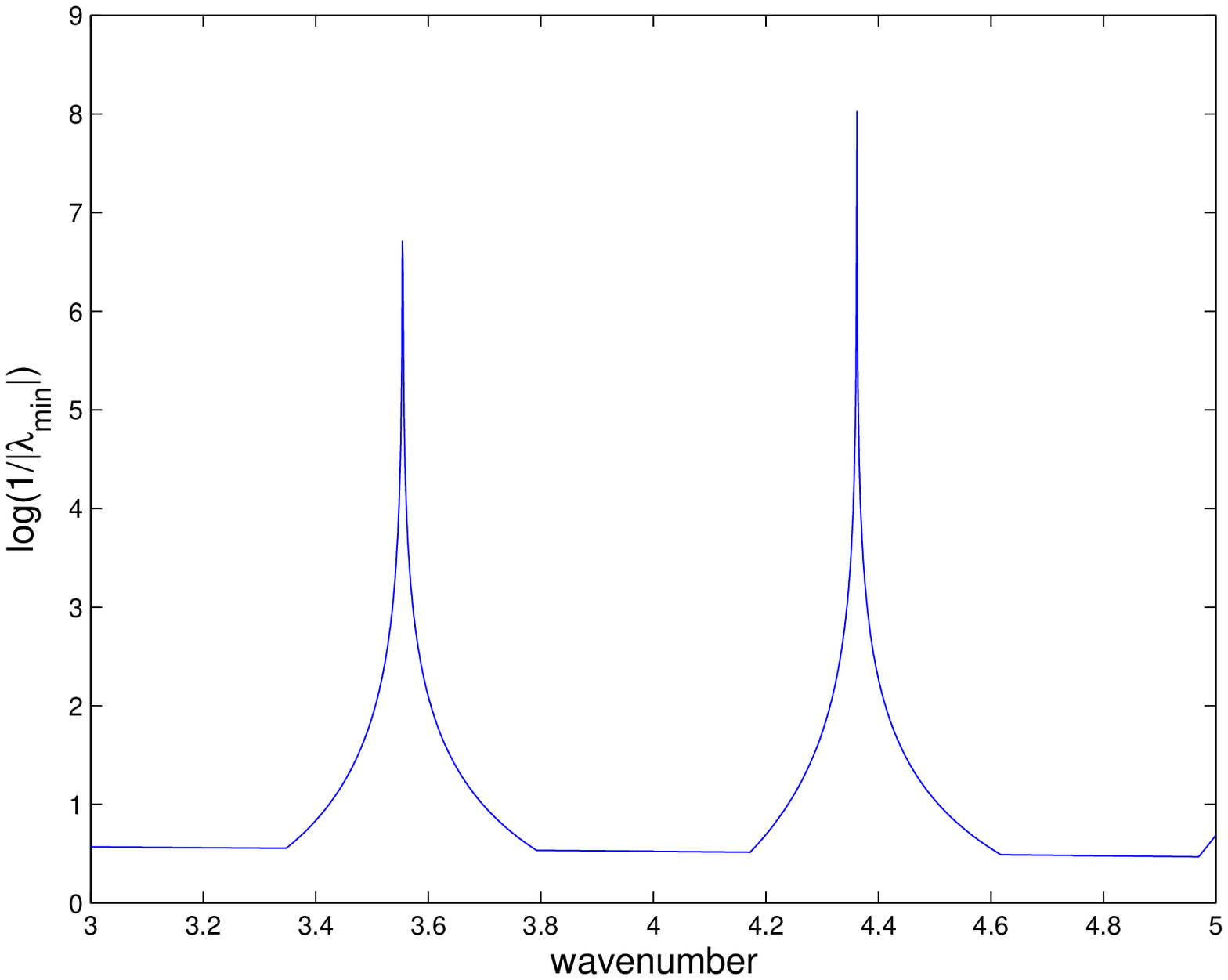}}
\end{tabular}
\end{center}
\caption{Log plot of $1/|\lambda_{min}|$. Left: $n=16$. Right: $n=9$.}
\label{haddarn}
\end{figure}

We also show the log plot of $1/|\lambda_{min}|$ by the method of  \cite{Cossoniere2011} in Fig.~\ref{haddarn}. Comparing 
Figures \ref{n16} and \ref{n9} with Figure \ref{haddarn}, it is clear that the probing method has much narrower span.

\section{Conclusions and Future Works}
In this paper, we proposed a probing method based on contour integrals for the transmission eigenvalue problem.
The method only tests if a given region contains an eigenvalue or not. Comparing with the existing methods, it needs little
a prior spectrum information and seems to be more efficient. 
The method can be viewed as an eigensolver without computing eigenvalues.
One advantage of the contour integral method is that it is suitable for parallel computing. Therefore, even the desired 
eigenvalues are dispersed, one can use a parallel scheme to capture them 
simultaneously. 

Note that one needs to construct two matrices for each wavenumber. It is time consuming if one wants to divide the
searching interval into more subintervals to improve accuracy. The work load is much more in three dimension.
Currently, we are developing a parallel version of the method using graphics processing units (GPUs).

\section*{Acknowlegement}
The work of F. Zeng is partially supported by the NSFC Grant (11501063). 
The work of J. Sun is supported in part by NSF DMS-1521555 and the US Army Research Laboratory and the US Army Research Office under the cooperative agreement number W911NF-11-2-0046. The work of L. Xu is partially supported by the NSFC Grant (11371385),  the Start-up fund of Youth 1000 plan of China and that of Youth 100 plan of Chongqing University.

\end{document}